\documentclass{article}
\usepackage{amssymb, amsthm, amsmath}
\usepackage{graphicx}
\usepackage{authblk}
\usepackage{float}
\usepackage{cleveref}

\newtheorem{theorem}{Theorem}
\newtheorem{lemma}{Lemma}
\newtheorem{proposition}{Proposition}

\begin{document}
\sloppy
\title{Classification of NMS-flows with unique twisted saddle orbit on orientable 4-manifolds} 
\author[1, 2]{Vladislav Galkin\footnote{ORCID: 0000-0001-6796-9228}}
\author[1, 2]{Olga Pochinka\footnote{ORCID: 0000-0002-6587-5305}}
\author[1, 2]{Danila Shubin\footnote{ORCID: 0000-0002-8495-4826; Corresponding author: schub.danil@yandex.ru}}
\date{}
\affil[1]{National Research University Higher School of Economics; 25/12, Bolshaya Pecherskaya St., Nizhny Novgorod, Russia, 603155}
\affil[2]{Saint Petersburg University, 7/9 Universitetskaya nab., St. Petersburg, 199034 Russia}
\maketitle
	
\begin{abstract} 
	Topological equivalence of Morse-Smale flows without fixed points (NMS-flows) under assumptions of different generalities was studied in a number of papers. 
	In some cases when the number of periodic orbits is small, it is possible to give exhaustive classification, namely to provide the list of all manifolds that admit flows of considered class, find complete invariant for topological equivalence and introduce each equivalence class with some representative flow. 
	This work continues the series of such articles. 
	We consider the class of NMS-flows with unique saddle orbit, under the assumption that it is twisted, on closed orientable 4-manifolds and prove that the only 4-manifold admitting the considered flows is the manifold $\mathbb S^3\times\mathbb S^1$. 
	Also, it is established that such flows are split into exactly eight equivalence classes and construction of a representative for each equivalence class is provided.
\end{abstract}
	
\section{Introduction and main results}

In the present paper we consider \emph{NMS-flows} $f^t$, namely \emph{non-singular} (without fixed points) Morse-Smale flows which are defined on orientable-manifold $M^4$. 
Non-wandering set of such flow consist of a finite number of hyperbolic periodic orbits. 
Asimov proved~\cite{Azimov} that ambient manifold of such flow is a union of round handles. However, if the number of orbits is small, topology of the ambient manifold can be specified. 
For instance, only lens spaces admit NMS-flows with two periodic orbits. 
Moreover, it was shown in~\cite{PoSh} that each lens space admit exactly two classes of topological equivalence except for 3-sphere $\mathbb S^3$ and projective space $\mathbb R\text{P}^3$ which both admit the unique equivalence class.

Campos et al.~\cite{CamposCorderoMartinezAlfaroVindel} argued that lens spaces are the only prime 3-manifolds that are ambient for NMS-flow with unique saddle periodic orbit, but this is not so. 
There exists infinite series of mapping tori non-homeomorphic to lens spaces which admit such flows~\cite{Shu21}. 
Moreover, necessary and sufficient conditions for topological equivalence of such flows were obtained in~\cite{PoSh-eqv}. 
Finally, any 3-manifold admitting such flows is a lens space or connected sum of two lens spaces or a small Seifert fibered 3-manifold~\cite{PoSh-top}. 
Note, that these results cannot be deduced from the classification of Morse-Smale flows with finite number of singular trajectories made by Umanskii~\cite{Umansky}. 

In the present paper we establish the topology of orientable 4-manifolds that are ambient for NMS-flows with exactly one saddle orbit assuming that it  is \emph{twisted} (its invariant manifolds are non-orientable). Besides, the number of equivalence classes of such flows on each admissible manifold is calculated.

Let us proceed to the formulation of the results.

Let $M^4$ be connected closed orientable 4-manifold, $f^t\colon M^4\to M^4$ be NMS-flow and $\mathcal O$ be its periodic orbit. 
There exists tubular neighborhood $V_{\mathcal O}$ homeomorphic to $\mathbb D^3\times \mathbb S^1$ such that the flow is topologically equivalent to the suspension over a linear diffeomorphism of the plane defined by the matrix which determinant is positive and eigenvalues are different from $\pm 1$ (see.~\Cref{Irw}). 
If absolute values of both eigenvalues are greater (less) than one, the corresponding periodic orbit is \emph{attracting (repelling)}, otherwise it is \emph{saddle}. 
The saddle orbit is called \emph{twisted} if both eigenvalues are negative and \emph{non-twisted} otherwise.

Consider the class $G^-_3(M^4)$ of NMS-flows $f^t\colon M^4\to M^4$ with unique saddle orbit which is twisted. 
Since the ambient manifold $M^4$ is the union of stable (unstable) manifolds of its periodic orbits, the flow $f^t\in G^-_3(M^4)$ has at least one attracting and at least one repelling orbit. 
In~\Cref{dyd} the following fact is established.
\begin{lemma}\label{lem:RAS} 
	The non-wandering set of any flow $f^t\in G^-_3(M^4)$ consists of exactly three periodic orbits $S,A,R$, saddle, attracting and repelling, respectively.
\end{lemma}

Saddle orbit $S$ of the flow $f^t\in G^-_3(M^4)$ can be either 3- or 2-dimensional; let $G^{-1}_3(M^4),\ G^{-2}_3(M^4)$ denote corresponding subclasses of $G^-_3(M^4)$. 
Obviously, since dimension of unstable manifold is invariant under an equivalence homeomorphism no flow in $G^{-1}_3(M^4)$ is topologically equivalent to any flow in $G^{-2}_3(M^4)$. 
Furthermore, $G^{-2}_3(M^4)=\{f^{-t}\colon f^t\in G^{-1}_3(M^4)\}$ and the flows $f^t,f'^t$ are topologically equivalent if and only if $f^{-t},f'^{-t}$ are topologically equivalent. 
This immediately implies that classification in the class $G^-_3(M^4)$ reduces to classification in the subclass $G^{-1}_3(M^4)$. 

Let $f^t\in G^{-1}_3(M^4)$. Since the flow $f^t$ in some tubular neighborhood of is topologically equivalent to the suspension over linear diffeomorphism of the plane, the topology of periodic orbits $A,S,R$  stable and unstable manifolds is:
\begin{itemize}
 	\item $W^u_S\cong \mathbb R^2\tilde\times\mathbb S^1$ (open solid Klein bottle);
 	\item $W^s_S\cong\mathbb R\tilde\times\mathbb S^1$ (open M\"{o}bius band);
 	\item $W^s_A\cong W^u_R\cong\mathbb R^2\times\mathbb S^1$ (open solid torus);
 	\item $W^u_A\cong W^s_R\cong\mathbb S^1$ (circle).
 \end{itemize} 

Let $\mathcal O\in\{S,A,R\}$. 
Choose the generator $\mathcal G_{\mathcal O}$ of boundary $T_{\mathcal O}=\partial V_{\mathcal O}\cong\mathbb S^2\times\mathbb S^1$ fundamental group which is homologous to $\mathcal O$ in $V_{\mathcal O}\cong\mathbb D^3\times\mathbb S^1$. 
By definition the manifold $T_S$ is secant for all flow trajectories except the periodic ones. Since the flow in some tubular neighborhood of $S$ is topologically equivalent to suspension, the set $K_S=W^u_S\cap T_S$ is homeomorphic to the Klein bottle.
Let $\lambda_{S},\mu_{S}$ be the knots (simple closed curves), which are generators of the fundamental group $\pi_1(K_S)$ with relation $[\lambda_{S}*\mu_{S}]=[\mu_{S}^{-1}*\lambda_{S}]$. 
We will call the curve $\mu_{S}$ \emph{meridian} and the curve $\lambda_{S}$ \emph{longitude}. 
By virtue of~\Cref{prop:uniq-Klein} the Klein bottle longitude embedded in $\mathbb S^2\times\mathbb S^1$ is a generator of fundamental group  $\pi_1(\mathbb S^2\times\mathbb S^1)$. 
Consider the longitude $\lambda_{S}$ be oriented in such way that its homotopy type $\langle\lambda_{S}\rangle$ in $T_S$ coincide with type $\langle \mathcal G{_S}\rangle$.
So, the set $K_A=W^u_S\cap T_A$ is the Klein bottle with longitude $\lambda_A$ which is pointwise transferred along the flow $f^t$ orbits from $\lambda_{S}$.

Since the flow in some tubular neighborhood of $S$ is topologically equivalent to suspension the set $\gamma_{S}=W^s_S\cap T_S$ is a knot in $T_S$, wrapping around $\mathcal G_S$ twice. 
We will assume that the knot $\gamma_{S}$ is oriented in such way that its homotopy type $\langle\gamma_{S}\rangle$ on $T_S$ coincides with homotopy type of $2\langle \mathcal G{_S}\rangle$. 
So the set $\gamma_{R}=W^s_S\cap T_R$ is a knot in $T_R$ and its orientation is induced by the flow $f^t$ from $\gamma_{S}$.
\begin{lemma}\label{lem:gen} Let $f^t\in G^{-1}_3(M^4)$ then the following conditions hold:
	\begin{enumerate}
		\item $\langle\lambda_{A}\rangle=\delta_A\langle \mathcal G_{A}\rangle,\,\delta_A\in\{-1,+1\}$ in $T_A$;
		\item $\langle\gamma_{R}\rangle=\delta_R\langle \mathcal G_{R}\rangle,\,\delta_R\in\{-1,+1\}$ in $T_R$.
	\end{enumerate}
\end{lemma}
Let 
$$C_{f^t}=(\delta_A,\delta_R).$$

\begin{theorem}\label{th:top-eqv}
	Flows $f^t,\, f'^t\in G^{-1}_3(M^4)$ are topologically equivalent if and only if $C_{f^t}=C_{f'^t}$.
\end{theorem}
\begin{theorem}\label{th:reali}
	For any element $C\in\mathbb S^0\times\mathbb S^0$ there exists a flow $f^t\in G_3^{-1}(M^4)$ such that $C=C_{f^t}$.
\end{theorem}
\begin{theorem}\label{th:topology}
	The only 4-manifold that is ambient for a flow of the class $G^-_3(M^4)$ is $\mathbb S^3\times\mathbb S^1$. Moreover $G^-_3(\mathbb S^3\times\mathbb S^1)$
	consists of eight classes of topological equivalence.
\end{theorem}
Note that weakening the saddle orbit twistedness condition fundamentally changes the picture. 
For example, in \cite{PoSh-wild} non-singular flows that are suspensions over Pixton diffeomorphisms on a three-dimensional sphere are considered. 
It is proved that in the class under consideration there exist flows with wildly embedded invariant manifolds of the saddle orbit. 
Moreover, there are an infinite number of topological equivalence classes for such flows.

\textit{Acknowledgements. } This work was performed at the Saint Petersburg Leonhard Euler International Mathematical Institute and supported by the Ministry of Science and Higher Education of the Russian Federation (agreement no. 075-15-2022-287).

\section{Flows of the class $G^{-1}_3(M^4)$}\label{dyd}
\subsection{Structure of periodic orbits}
This section is devoted to proof of~\Cref{lem:RAS}: non-wandering set of any flow $f^t\in G^{-1}_3(M^4)$ consists of three periodic orbits $S,A,R$, saddle, attracting and repelling respectively.
\begin{proof}
The proof is based on the following representation of the ambient manifold $M^4$ of the NMS-flow $f^t$ with the set of periodic orbits $Per_{f^t}$ (see, for example, \cite{Sm})

\begin{equation}\label{Mob}
	M^4 = \bigcup\limits_{\mathcal O \in Per_{f^t}} W^u_{\mathcal O}=\bigcup\limits_{\mathcal O \in Per_{f^t}} W^s_{\mathcal O},
\end{equation} as well as the asymptotic behavior of invariant manifolds 
\begin{equation}\label{neust}
	{\rm cl}(W^u_{\mathcal O}) \setminus W^u_{\mathcal O} = \bigcup\limits_{\tilde{\mathcal O} \in Per_{f^t}\colon W^u_{\mathcal O}\cap W^s_{\mathcal O}\neq \varnothing} W^u_{\tilde{\mathcal O}},
\end{equation}
\begin{equation}\label{ust}
	{\rm cl}(W^s_{\mathcal O}) \setminus W^s_{\mathcal O} = \bigcup\limits_{\tilde{\mathcal O} \in Per_{f^t}\colon W^s_{\mathcal O}\cap W^u_{\mathcal O}\neq \varnothing} W^s_{\tilde{\mathcal O}}.
\end{equation}
In particular, it follows from \cref{Mob} that any NMS-flow has at least one attracting orbit and at least one repulsive one. Moreover, if an NMS-flow has a saddle periodic orbit, then the basin of any attracting orbit has a non-empty intersection with an unstable manifold of at least one saddle orbit (see~Proposition~2.1.3~\cite{begin}) and a similar situation with the basin of a repulsive orbits.

Now let $f^t\in G_3^{-1}(M^3)$ and $S$ be its only saddle orbit. It follows from the relation (\ref{Mob}) that $W^u_S\setminus S$ intersects only basins of attracting orbits. Since the set $W^u_S\setminus S$ is connected and the basins of attracting orbits are open, then $W^u_S$ intersects exactly one such basin. Denote by $A$ the corresponding attracting orbit. Since there is only one saddle orbit, there is only one attracting orbit. Similar reasoning for $W^s_S$ leads to the existence of a unique repulsive orbit $R$.
\end{proof}

\subsection{Canonical neighborhoods of periodic orbits}\label{subsec:canon-hood}

Recall the definition of a superstructure. Let $\varphi\colon M^3\to M^3$ be a diffeomorphism of a 3-manifold. We define the diffeomorphism $g_{\varphi}\colon M^3\times \mathbb R^1 \to M^3\times \mathbb R^1$ by the formula $$g_{\varphi}(x_1,x_2,x_3, x_4) = (\varphi(x_1,x_2, x_3),x_4-1).$$
Then the group $\{g_{\varphi}^n\}\cong\mathbb Z$ acts freely and discontinuously on $M^3\times \mathbb R^1$, whence the orbit space $\Pi_\varphi = M ^3\times \mathbb R^1 g_{\varphi}$ is a smooth 4-manifold, and the natural projection $v_\varphi\colon M^3\times \mathbb R^1\to \Pi_\varphi$ is a covering.
At the same time, the flow $\xi^t\colon M^3\times \mathbb R^1\to M^3\times \mathbb R^1$, given by the formula $$\xi^t(x_1, x_2, x_3, x_4)=(x_1, x_2, x_3, x_4+t),$$ induces the flow $[\varphi]^t= v_\varphi \xi^t v^{-1}_\varphi\colon\Pi_\varphi\to \Pi_\varphi$. The flow $[\varphi]^t$ is called the \emph{suspension of the diffeomorphism $\varphi$}.

We define the diffeomorphisms $a_0,a_1,a_2,a_3\colon \mathbb R^{3}\to \mathbb R^{3}$ by the formulas
$$a_3(x_1, x_2, x_3) = (2x_1, 2x_2, 2x_3),\ a_0 = a_3^{-1},$$
$$a_{\pm 1}(x_1, x_2, x_3) = (\pm 2x_1, \pm 1/2x_2, 1/2 x_3),\ a_{\pm 2} = a_{\pm 1}^{- 1}.$$
Let
\begin{gather*}
	V_{0}= \{ (x_1,x_2, x_3, x_4)\in \mathbb R^3 |\ 4^{x_4} x_1^2 + 4^{x_4}x^2_2 + 4^{x_4}x^ 2_3 \leqslant 1 \},\\
	V_{\pm 1}= \{ (x_1,x_2, x_3, x_4)\in \mathbb R^3 |\ 4^{-x_4} x_1^2 + 4^{x_4}x^2_2 + 4^{x_4 }x^2_3 \leqslant 1 \},\\
	V_{\pm 2}= \{ (x_1,x_2, x_3, x_4)\in \mathbb R^3 |\ 4^{-x_4} x_1^2 + 4^{-x_4}x^2_2 + 4^{ x_4}x^2_3 \leqslant 1 \},\\
	V_{3}= \{ (x_1,x_2, x_3, x_4)\in \mathbb R^3 |\ 4^{-x_4} x_1^2 + 4^{-x_4}x^2_2 + 4^{-x_4 }x^2_3 \leqslant 1 \}.
\end{gather*}
For $i\in\{0,\pm 1,\pm 2,3\}$ we set $v_i = v_{a_i}$, $T_i=\partial V_i$ and $ \mathbb V_i=v_i(V_i),\ \mathbb T_{i} = v_i(T_i)$.

The following statement, proved by M.~Irwin~\cite{Irwin}, describes the behavior of flows in a neighborhood of hyperbolic periodic orbits.

\begin{proposition}[M.~Irwin~\cite{Irwin}]\label{prop:irwin-eqv} If $\mathcal O$ is a hyperbolic orbit of a flow $f^t\colon M^4\to M^4$ defined on an orientable 4-manifold $M ^4$, then there exists a tubular neighborhood $V_{\mathcal O}$ of the orbit $\mathcal O$ such that the flow $f^t\big|_{V_{\mathcal O}}$ is topologically equivalent, by means of some homeomorphism $ H_{\mathcal O}$, to one of the following streams:
	\begin{itemize}
		\item $[a_0]^t|_{\mathbb V_0}$ if $\mathcal O$ is an attracting orbit;
		\item $[a_{ 1}]^t|_{\mathbb V_{1}}$ if $\mathcal O$ is a non-twisted saddle orbit with a two-dimensional unstable manifold $W^u_{\mathcal O}$;
		\item $[a_{-1}]^t|_{\mathbb V_{-1}}$ if $\mathcal O$ is a twisted saddle orbit with a two-dimensional unstable manifold $W^u_{\mathcal O}$ ;
		\item $[a_{ 2}]^t|_{\mathbb V_{2}}$ if $\mathcal O$ is a non-twisted saddle orbit with an unstable 3-manifold $W^u_{\mathcal O}$;
		\item $[a_{-2}]^t|_{\mathbb V_{-2}}$ if $\mathcal O$ is a twisted saddle orbit with an unstable 3-manifold $W^u_{\mathcal O}$ ;
		\item $[a_3]^t|_{\mathbb V_3}$ if $\mathcal O$ is a repelling orbit.
	\end{itemize}
	\label{Irw}
\end{proposition}

The neighborhood $V_\mathcal O=H_{\mathcal O}(\mathbb V_i),\,i\in\{0,\pm 1, \pm2, 3\}$ described in~\Cref{Irw} is called \emph{canonical neighborhood} of the periodic orbit $\mathcal O$.

When proving topological equivalence, we will use the following fact, which follows from the proof of Theorem~4 and Lemma~4 in~\cite{PoSh}, and can also be found in~\cite{Umansky}~(Theorem~1.1).

\begin{proposition}\label{prop:h->H} A homeomorphism $h_i\colon \partial \mathbb V_i\to \partial \mathbb V_i$ for $i\in\{0,3\}$ extends to a homeomorphism $H_i\colon \mathbb V_i\to \mathbb V_i$ realizing the equivalence flow $[a_i]^t$ with itself if and only if the induced isomorphism $h_{i*}\colon \pi_1(\partial \mathbb V_i)\to \pi_1(\partial \mathbb V_i)$ is identical.
\end{proposition}

\subsection{Trajectory mappings}\label{potr}
Consider a flow $f^t\colon M^4\to M^4$ from the set $G^{-1}_3(M^3)$. Then $V_A=H_A(\mathbb V_0),\,V_R=H_R(\mathbb V_3),\,V_S=H_S(\mathbb V_{-2})$. 
Let $\Gamma = \{(x_1, x_2, x_3, x_4)\in T_{-2}|\ 4^{x_4}x_3^2 = 1/2\},\ \Gamma^u = Ox_2x_3x_4\cap T_ {-2},\ \Gamma^s = Ox_1x_4\cap T_{-2}$. 
By construction, the set $T_{-2}$ is homeomorphic to $\mathbb S^2\times\mathbb R$, the set $\Gamma$ consists of two surfaces, each of which is homeomorphic to $\mathbb S^1\times\mathbb R$ , dividing $T_{-2}$ into three connected components, one of which $N^u$ contains the cylinder $\Gamma^u\cong\mathbb S^1\times\mathbb R$, and the union $N^s$ the other two contains a pair of $\Gamma^s\cong\mathbb S^0\times\mathbb R$ curves, one curve in each component. Then on $T_S$
\begin{itemize}
	\item $K_S=H_S(v_{-2}(\Gamma^u)))$ is a Klein bottle;
	\item $\gamma_{S}=H_S(v_{-2}(\Gamma^s)))$ is a knot winding twice around the generator $\mathcal G_S$;
	\item $N^u_S = H_S(v_{-2}({\rm cl}(N^u)))$ is a tubular neighborhood of $K_S$;
	\item $N^s_S = H_S(v_{-2}({\rm cl}(N^s)))$ is a tubular neighborhood of $\gamma_S$;
	\item $\partial N^u_S=\partial N^s_S= H_S(v_{-2}((\Gamma)))$ is a two-dimensional torus.
\end{itemize}
Let $$N^s_R=\left(\bigcup\limits_{t>0,\,w\in N^s_S}f^{-t}(w)\right)\cap T_R,\,\,N^ u_R=T_R\setminus N^s_R,$$ $$N^u_A=\left(\bigcup\limits_{t>0,\,w\in N^u_S}f^{t}(w)\right)\cap T_A,\,\,N^s_A=T_A\setminus N^u_A$$ and introduce the following mappings:
\begin{itemize}
	\item define a continuous function $\tau_{R}\colon T_R\to\mathbb R^+$ so that $f^{\tau_{R}(r)}(r)\in N^s_S$ for $r\ in N^s_R$ and the set $\mathcal T=\bigcup\limits_{r\in N^u_R}f^{\tau_{{R}}(r)}(r)$ does not intersect the torus $T_A$, we set ${\mathcal T}_R={\mathcal T}\cup N^s_S$ and define a homeomorphism $\psi_{R}\colon T_R\to{\mathcal T}_R$ by the formula
	$\psi_{R}(r)=f^{\tau_{{R}}(r)}(r)$, denote by $\mathcal V_R$ the closure of the connected component of the set $M^4\setminus \mathcal T_R$ , containing $R$;
	
	\item we set ${\mathcal T}_A={\mathcal T}\cup N^u_S$, define a continuous function $\tau_{A}:T_A\to\mathbb R^+$ so that $f^{- \tau_{A}(a)}(a)\in {\mathcal T}_A$ for $a\in T_A$ and define a homeomorphism $\psi_{A}:T_A\to{\mathcal T}_A$ by the formula
	$\psi_{A}(a)=f^{-\tau_{{A}}(a)}(a)$, denote by $\mathcal V_A$ the closure of the connected component of the set $M^4\setminus \mathcal T_A $ containing $A$.
\end{itemize}

We will call the introduced homeomorphisms $\psi_{R},\ \psi_{A}$ \emph{trajectory maps}. Note that the ambient manifold $M^4$ is represented as
$$M^4 = \mathcal V_A \cup V_S \cup \mathcal V_R.$$
Moreover, in the manifolds $\mathcal V_A$ and $\mathcal V_R$ the flow $f^t$ is topologically equivalent to the suspensions $[a_0]^t$ and $[a_3]^t$, respectively.

\section{Homotopy types of knots $\lambda_{A},\gamma_{R}$}

In this section, we will prove~\Cref{lem:gen}. To do this, we first describe the properties of the embedding of the Klein bottle into the manifold $\mathbb S^2\times\mathbb S^1$.

Recall that the Klein bottle $\mathbb K$ is the square $[0, 1]\times [0, 1]$ with sides glued by the relation
$$(x, 0)\sim (x, 1),\quad (0, y)\sim (0, 1-y).$$ 
Let $v\colon [0, 1]\times [0, 1] \to \mathbb K$ be the natural projection, then the curves $$\lambda = v([0, 1]\times\{1/2 \}),\,\mu = v(\{0\}\times[0, 1])$$ are generators of the fundamental group $\pi_1(\mathbb K)$ with relation $$[\lambda*\mu] =[\mu^{-1}*\lambda],$$ where the curve $\lambda$ is called \emph{longitude} and the curve $\mu$ is called \emph{meridian}.

It is well known that the Klein bottle does not embed into $\mathbb R^3$, however, it can be embeddable into $\mathbb S^2\times \mathbb S^1$, for example by defining the embedding $\tilde e_0\colon [0 , 1]\times [0, 1]\to \mathbb S^2\times\mathbb S^1$ by the formula
$$\tilde e_0(x, y) = \left(\sin \pi x \cos 2 \pi y,\ \cos \pi x \cos 2 \pi y,\ \sin 2 \pi y,e^{ 2\pi i x}\right)$$ and noticing that $\tilde e_0(x,y)=\tilde e_0(x',y')\iff(x,y)\sim(x',y' )$. Then (see,~for~example,~\cite[Chapter 5]{Kosn}) $$e_0=\tilde e_0v^{-1}\colon\mathbb K\to\mathbb S^2\times\mathbb S^1$$  is the desired embedding of the Klein bottle in $\mathbb S^2\times\mathbb S^1$. Let 
$$K_0 =e_0(\mathbb K).$$
\begin{figure}[h!]
\centerline{\includegraphics[width=0.5\textwidth]{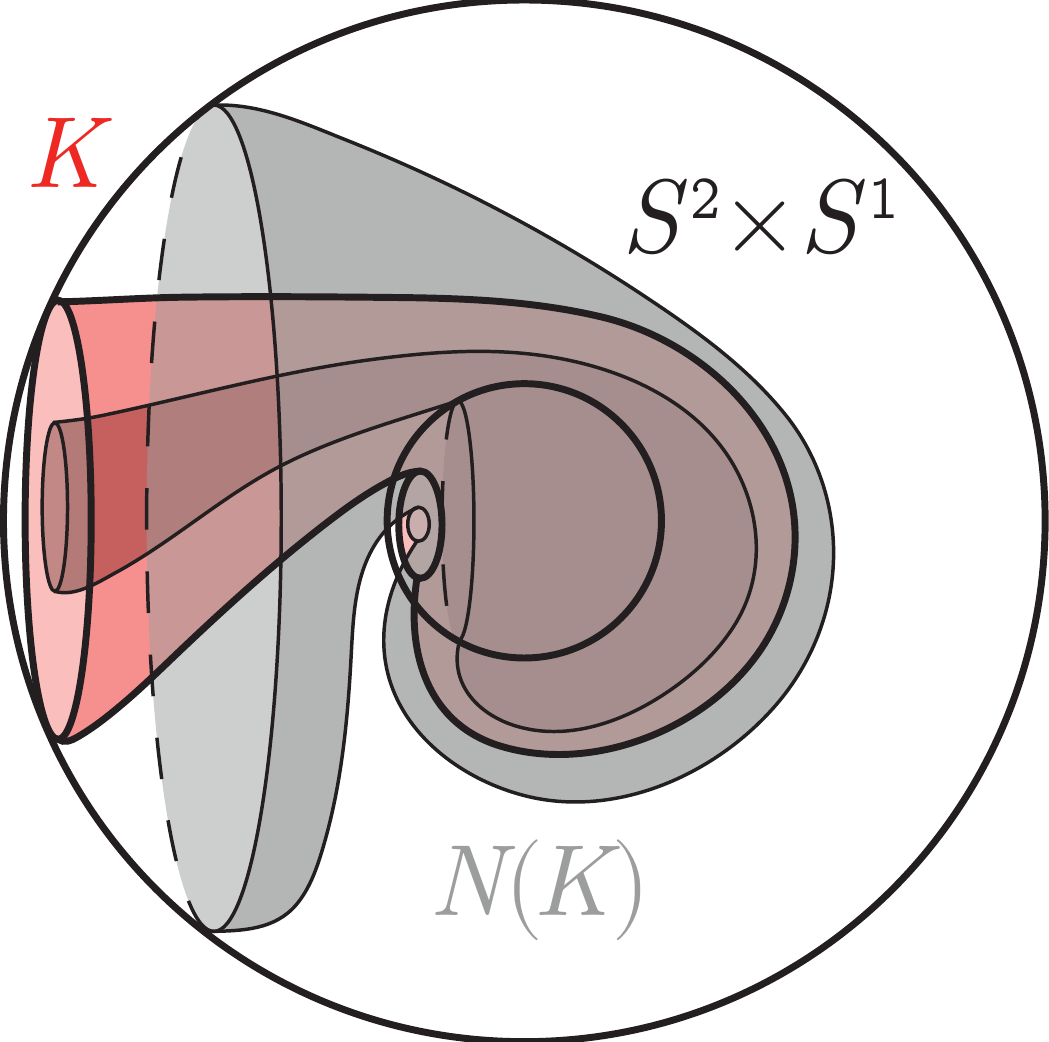}}
	\caption{Bottle of Klein in $\mathbb S^2\times\mathbb S^1$}\label{pic:Klein}
\end{figure}
\begin{proposition}[Proposition 1.4, \cite{twoLinks}]\label{prop:uniq-Klein} Let $e\colon \mathbb K\to\mathbb S^2\times \mathbb S^1$ be an embedding Klein bottles $\mathbb K$, $K=e(\mathbb K)$, $N(K)\subset\mathbb S^2\times \mathbb S^1$ be a tubular neighborhood of $K$ and $V( K)=\mathbb S^2\times \mathbb S^1\setminus {\rm int}\,N(K)$ (see~\Cref{pic:Klein}). Then:
	\begin{itemize}
		\item[1)] the curve $e(\lambda)$ is a generator of the fundamental group $\pi_1(\mathbb S^2\times \mathbb S^1)$;
		\item[2)] the set $V(K)$ is a solid torus whose meridian is homotopic to the curve $e(\mu)$;
		\item[3)] there exists an orientation-preserving homeomorphism $h\colon \mathbb S^2\times \mathbb S^1\to \mathbb S^2\times \mathbb S^1$ such that $h(K) = K_0$ and $h_* = id\colon \pi_1(\mathbb S^2\times \mathbb S^1)\to \pi_1(\mathbb S^2\times \mathbb S^1)$.
	\end{itemize}
\end{proposition}

It remains to prove~\Cref{lem:gen}. To do this, recall that we have chosen $T_{\mathcal O}=\partial V_{\mathcal O}\cong\mathbb S^2\times\mathbb S^1,\,\mathcal O\in\{A ,S,R\}$ generator $\mathcal G_{\mathcal O}$ of the fundamental group of $T_{\mathcal O}$, homologous in $V_{\mathcal O}\cong\mathbb D^3\times\mathbb S^1$ orbit $\mathcal O$. 
Due to the fact that canonical neighborhoods of periodic orbits can be chosen so that $M^4 = \mathcal V_A \cup V_S \cup \mathcal V_R$ (see~\Cref{potr}), everywhere below we assume that $ V_{A} = \mathcal V_{A},\ V_{R} = \mathcal V_{R}$.

We also established that the set $K_S=W^u_S\cap T_S$ is a Klein bottle on $T_S\cong\mathbb S^2\times\mathbb S^1$ and oriented its parallel $\lambda_{S}$ so that 
$\langle\lambda_{S}\rangle=\langle \mathcal G_{S}\rangle$ on $T_S$. Since the set $K_A=W^u_S\cap T_A$ coincides with $K_S$, then $\lambda_{A}=\lambda_{S}$. 
We also established that the set $\gamma_{S}=W^s_S\cap T_S$ is a knot on $T_S$ and oriented so that $\langle\gamma_{S}\rangle=2\langle \mathcal G_{ S}\rangle$ to $T_S$. Since the set $\gamma_{R}=W^s_S\cap T_R$ coincides with $\gamma_{S}$, then $\gamma_{R}=\gamma_{S}$.

Let us show that the knots $\lambda_A, \gamma_R$ are generators in the fundamental groups of the manifolds $T_A, T_R$, respectively.

\begin{proof} Since $T_A$ is homeomorphic to the manifold $\mathbb S^2\times\mathbb S^1$ and $\lambda_A$ is a parallel of the Klein bottle $K_A\subset T_A$, then item 1) of~\Cref {prop:uniq-Klein} implies that $\lambda_A$ is a generator in the fundamental group of $T_A$, that is, on $T_A$
$$\langle\lambda_{A}\rangle=\delta_A\langle \mathcal G_{A}\rangle$$
for $\delta_A\in\{-1,+1\}$.
The set $N(K_A)=N^u_S$ is a tubular neighborhood of the Klein bottle $K_A$ in $T_A$. It follows from item 2) of~\Cref{prop:uniq-Klein} that the set $V(K_A)=T_A\setminus {\rm int}\,N(K_A)$ is a solid torus. The set $N(\gamma_{R})=N^s_S$ is a tubular neighborhood of the knot $\gamma_{R}$ in $T_R$. On the other hand $T_R\setminus {\rm int}\,N(\gamma_{R})=V(K_A)$. Thus, the complement to the tubular neighborhood of $\gamma_{R}$ in $T_R$ is a solid torus. By \cite[Proposition~4.2]{begin}, $\gamma_R$ is a generator in $T_R$ and, therefore,
$$\langle\gamma_{R}\rangle=\delta_R\langle \mathcal G_{R}\rangle$$
for $\delta_R\in\{-1,+1\}$.
\end{proof}

\section{Classification of flows of the set $G^{-1}_3(M^4)$}
In this section, we will prove Theorem~\ref{th:top-eqv}.

\begin{proof} $ $
	
\textit{Necessity.} Let flows $f^t$ and $f'^t$ have invariants $C_{f^t} = (\delta_A, \delta_R),\ C_{f'^t} = (\delta_ {A'}, \delta_{R'})$ and are topologically equivalent by the homeomorphism $H\colon M^4\to M^4$. Let us show that $C_{f^t}=C_{f'^t}$.

Let $h_A = H\big|_{T_A}$ and $T_{A'}=h(T_A)$. Then by~\Cref{prop:h->H}
$$\langle \mathcal G_{A'} \rangle =h_{A*}\langle \mathcal G_A\rangle.$$ It follows from item 3) of~\Cref{prop:uniq-Klein} that $\lambda_ {{A'}}=h_A(\lambda_{A})$ is the parallel of the Klein bottle $K_{A'}$. Since the longitude $\lambda_{A}$ of the Klein bottle is oriented consistent with the saddle orbit $S$ and $H$ transforms the orbit $S$ into the orbit $S'$ with orientation preservation, then $\lambda_{{A'}}$ is oriented consistent with the saddle orbit $S'$. On the other side,
by~\Cref{lem:gen},
$$\langle \lambda_A\rangle = \delta_{A} \langle \mathcal G_{A} \rangle,\,\langle \lambda_{A'} \rangle = \delta_{A'} \langle \mathcal G_{ A'} \rangle,$$
whence, by virtue of a simple chain of equalities
$$\delta_{A'} \langle\mathcal G_{A'}\rangle = h_{A*}( \delta_A \langle\mathcal{G}_A\rangle) =\delta_A \langle\mathcal G_{A' }\rangle, $$
we get that $\delta_A = \delta_{A'}$. It is proved similarly that $\delta_R = \delta_{R'}$. Thus $C_{f^t} = C_{f'^t}$.

\textit{Sufficiency.} Let the flows $f^t$ and $f'^t$ have equal invariants $C_{f^t} = (\delta_A, \delta_R),\ C_{f'^t} = (\delta_{A'}, \delta_{R'})$. Let us show that the flows $f^t$ and $f'^t$
are topologically equivalent.

\Cref{prop:irwin-eqv} implies that the homeomorphism $$H\big|_{V_S} = H^{}_{S'}H_S^{-1}\colon V_S\to V_{S' }$$ is topological equivalence homeomorphism of the flows $f^t\big|_{V_S}$ and $f'^t\big|_{V_{S'}}$. It remains to extend this homeomorphism to $\mathcal V_A$ and $\mathcal V_R$.

The homeomorphism $H$ is already defined on the set $\mathcal T_A\cap T_S$, which is a tubular neighborhood $N(K_A)$ of the Klein bottle $K_A$. By item 2) of~\Cref{prop:uniq-Klein} the set $V(K_A)=T_A\setminus {\rm int}\,N(K_A)$ is a solid torus whose meridian is homotopic to the meridian $\mu_A$ of the  Klein bottle $K_A$ to $N(K_A)$. It follows from the properties of the homeomorphism $H$ that
$K_{A'}=H(K_A)$ and $N(K_{A'})=H(N(K_A))$ is a tubular neighborhood of the Klein bottle $K_{A'}$. By point 2) of~\Cref{prop:uniq-Klein} the set $V(K_{A'})=T_{A'}\setminus {\rm int}\,N(K_{A'})$ is a solid torus whose meridian is homotopic to the meridian $\mu_{A'}$ of the Klein bottle $K_{A'}$ in $N(K_{A'})$. Since any homeomorphism of the Klein bottle does not change the homotopy class of the meridian (see,~for~example,~\cite[Lemma~5]{lickorish}), the homeomorphism $H\colon\partial V(K_A)\to\partial V(K_{A'} )$ extends to the homeomorphism $H\colon V(K_A)\to V(K_{A'})$ (see,~for example,~\cite[Exercise~2E5]{Rolfsen}).
Thus $H$ is defined on $\mathcal T_A$ and $\mathcal T_R$.

Since the parallel $\lambda_{A}\,(\lambda_{{A'}})$ of the Klein bottle is oriented consistent with the saddle orbit $S\,(S')$ and $H$ transforms the orbit $S$ into the orbit $S'$ orientation-preserving, then
$H_*(\langle \lambda_{A}\rangle) = \langle \lambda_{{A'}} \rangle$. Since $\delta_A=\delta_{A'}$, then $H_*(\delta_A\langle \lambda_{{A}}\rangle) =\delta_{A'} \langle \lambda_{{A'}} \rangle$ and hence $H_*(\langle \mathcal G_A\rangle) = \langle \mathcal G_{A'} \rangle$. 
By~\Cref{prop:h->H} the homeomorphism $H\big|_{\mathcal T_A}$ extends to $\mathcal V_A$ by  a homeomorphism realizing the equivalence of flows $f^t\big|_{\mathcal V_A}$ and $f'^t\big|_{\mathcal V_{A'}}$. 
Similarly, $H$ can be extended to $\mathcal V_R$. Thus, the homeomorphism $H$ is defined on the whole $M^4$ and realizes the equivalence of the flows $f^t,\,f'^t$.
\end{proof}

\section{Realization of flows by admissible set}

In this section, we will prove~\Cref{th:reali}: for any element $C\in\mathbb S^0\times\mathbb S^0$ there is a flow $f^t\in G_3^{-1}( M^4)$ such that $C=C_{f^t}$.

\begin{proof}
Let us construct the flow $f^t\colon \mathbb S^2\times \mathbb S^1\to \mathbb S^2\times \mathbb S^1$ with the invariant $C_{f^t} = (+1, + 1)$ as a suspension over the sphere diffeomorphism $\zeta\colon \mathbb S^2\to \mathbb S^2$ with three periodic orbits. To do this, we describe the construction of the diffeomorphism $\zeta$.

Let the function $\sigma\colon \mathbb R_+\to [0,1]$ be given by the formula (see~Fig.~\ref{pic:sigma}):
$$\sigma(x) = \begin{cases}
	\exp{\left(2 - \dfrac{1}{\left( \dfrac{1}{x} - 2\right)^2} - \dfrac{1}{(x - 2)^2} \right )}, & x\in\left[ \frac{1}{2}, 2\right] \\ 0, & x\notin \left[ \frac{1}{2}, 2\right].
\end{cases}$$
Note that $\sigma(1/x) = \sigma(x).$
\begin{figure}[h!]
	\centerline{\includegraphics[width=0.7\textwidth]{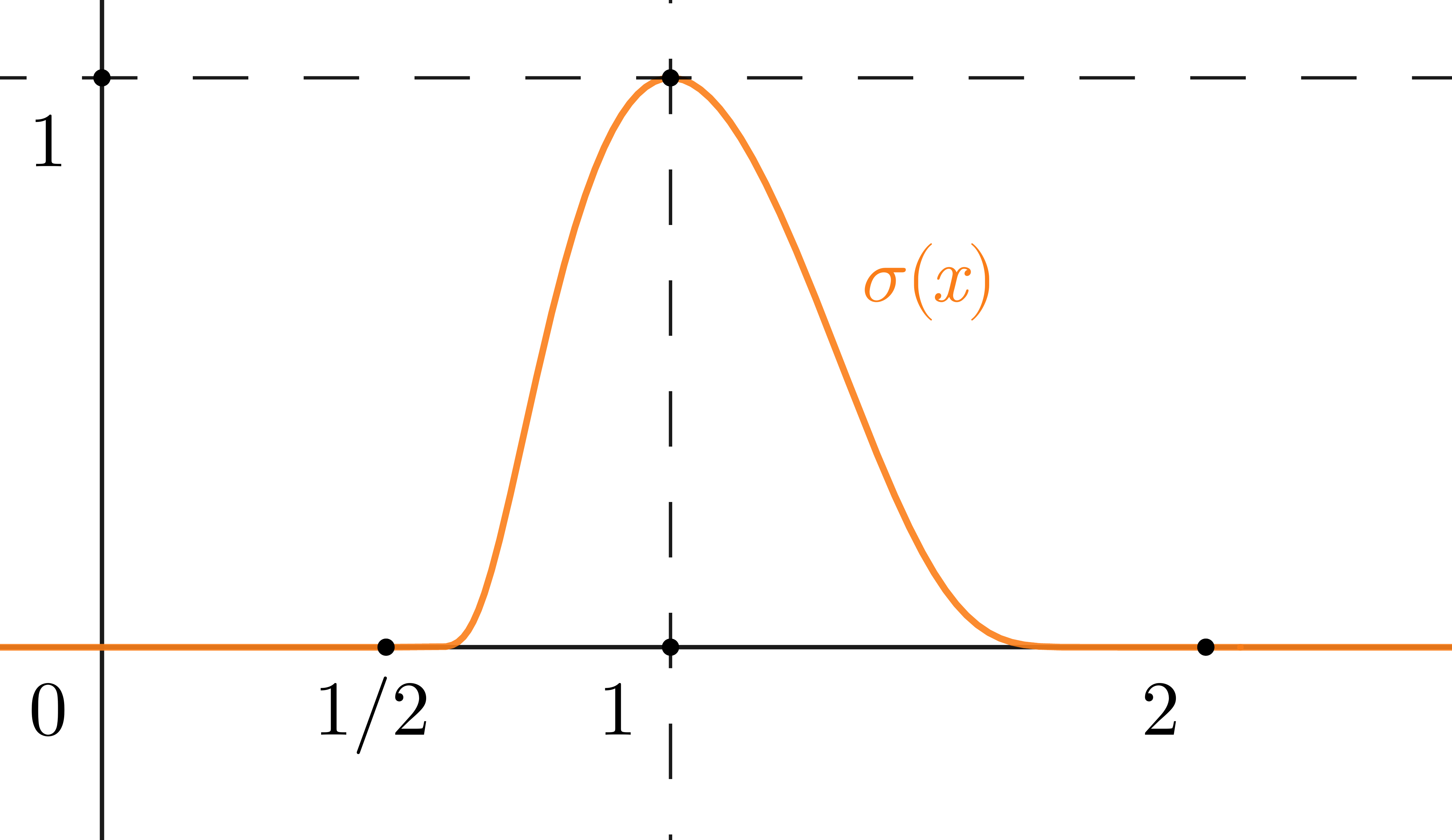}}
	\caption{Graph of the function $\sigma(x)$.}\label{pic:sigma}
\end{figure}

Introduce spherical coordinates in $\mathbb R^3$ by $x = r\cos\varphi\cos\theta,\ y = r\cos\varphi\sin\theta, z = r\sin\varphi,\quad r >0,\ \varphi\in[-\pi; \pi),\ \theta\in [-\pi/2; \pi/2)$.
Let $\chi^t\colon \mathbb R^3\to \mathbb R^3$ be the flow defined by system of equations:
$$\begin{cases}
	\dot r = \begin{cases}
		-r(r-1), & 0\leqslant r \leqslant 1\\
		1-r, & 1 \leqslant r
	\end{cases}\\
	\dot\varphi = 
	-\sigma(r)\sin \varphi\\
	\dot \theta = -\sigma(r)\theta
\end{cases}$$
and the diffeomorphism $q\colon \mathbb R^3\to \mathbb R^3$ be defined by the formula:
$$q(x, y, z) = (1/x,\ -1/y,\ 1/z).$$
\begin{figure}[h!] \centerline{\includegraphics[width=0.5\textwidth]{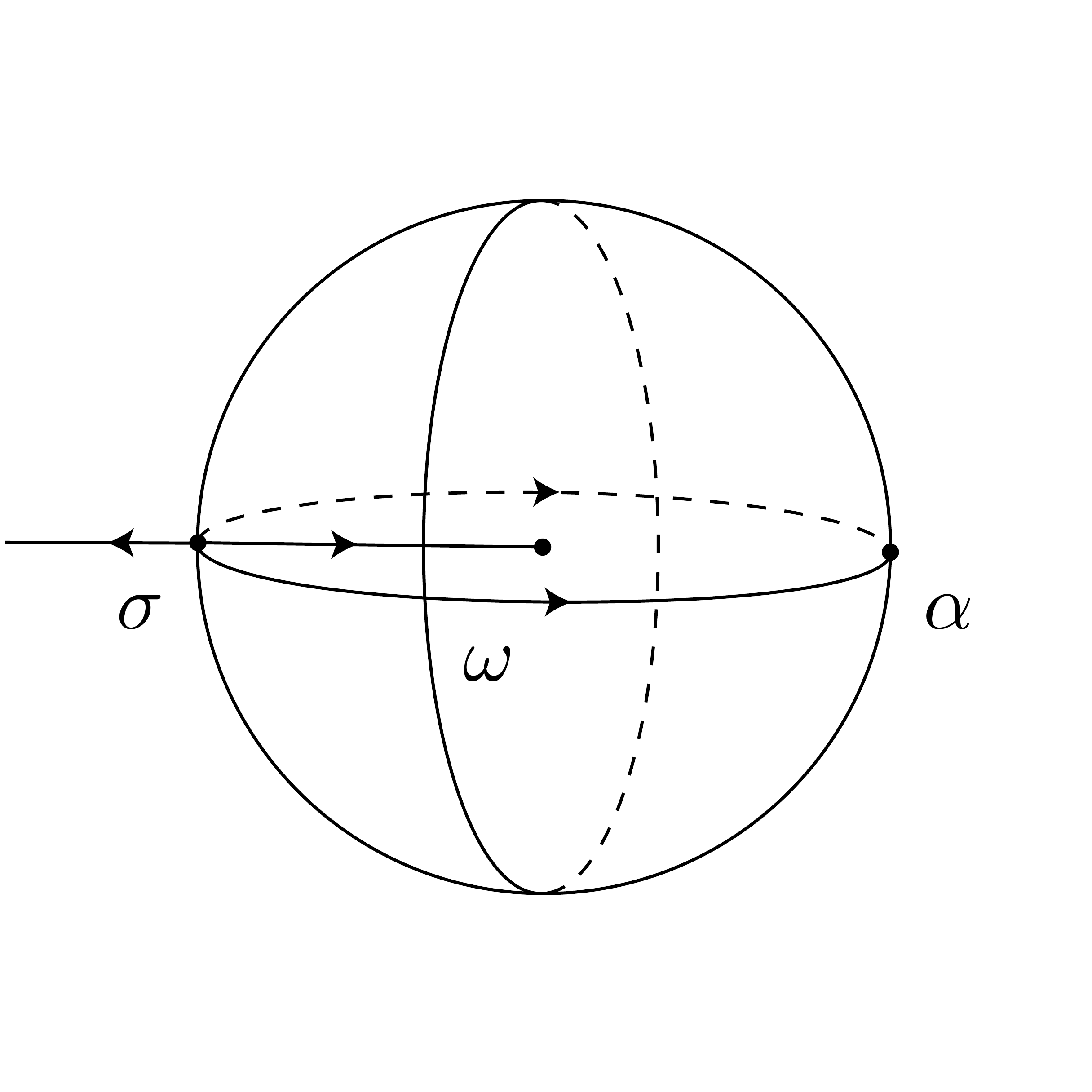}}		
	\caption{Flow $\chi^t$ phase portrait}\label{pic:phase-portrait}	
\end{figure}

Using stereographic projection (see~\Cref{pic:stereo}) $\vartheta\colon \mathbb S^3\setminus\{N\}\to \mathbb R^3$ ($N = (0 , 0, 0, 1), S = (0, 0, 0, -1)$) by the given formula:
\begin{equation*}\label{stereo+}
	\vartheta(x_1,x_2, x_3,x_4)=\left(\frac{x_1}{1-x_{4}}, \frac{x_2}{1-x_{4}}, \frac{x_3}{1 -x_{4}}\right).
\end{equation*}
\begin{figure}[h!]
\centerline{\includegraphics[width=1\textwidth]{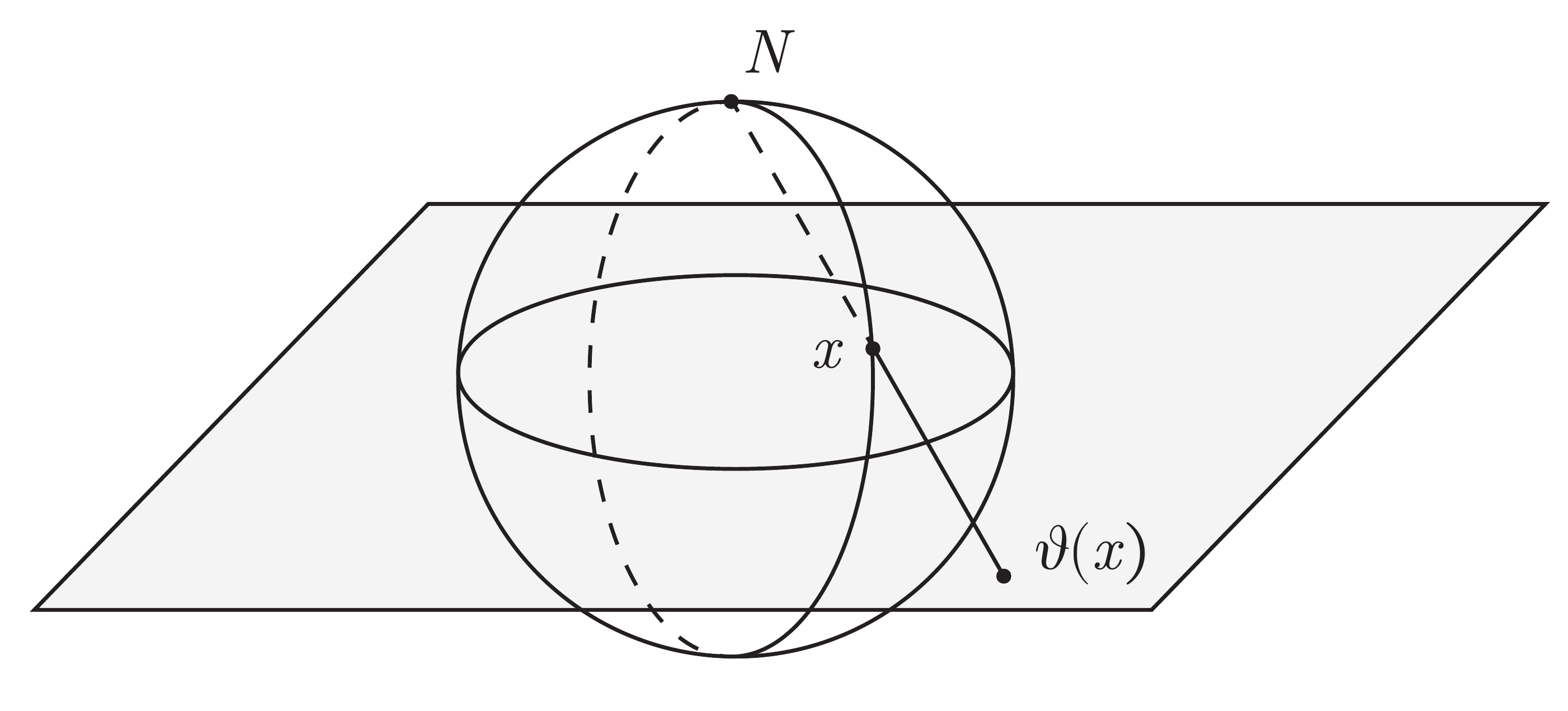}}
	\caption{Stereographic projection}\label{pic:stereo}
\end{figure}
project the diffeomorphism $q\chi^1$ onto $\mathbb S^3$:
$$f(x) = \begin{cases}
	\vartheta^{-1} q\chi^1 \vartheta(x), & x\notin \{N, S\}\\
	S, & x = N,\\
	N,&x=S
\end{cases}$$

Then the flow $f^t=[f]^t$ belongs to the class $G^{-1}_3(M^3)$ and $C_{f^t} = (+1, +1)$.
By construction, $f$ is an orientation-preserving diffeomorphism of the 3-sphere, and hence the ambient manifold of the suspension $[f]^t$ is homeomorphic to $\mathbb S^3\times \mathbb S^1$.

We construct the rest of the flows of class $G^{-1}_3(M^4)$ by modifying the constructed flow $f^t$.

Let $\mathcal O$ be an attractive or repelling periodic orbit of $f^t$ and $V_{\mathcal O} = V^0_{\mathcal O}$ be its canonical neighborhood, $V^{t }_{\mathcal O} =f^t(V^0_{\mathcal O})$. Without loss of generality, we assume that $V^{-1}_A\cap V^1_R=\emptyset$. Let $\vec v(x)$ denote the vector field induced by the flow $f^t$ on $\mathbb S^3\times \mathbb S^1$. For points $x$ that belong to the basin of the orbit $\mathcal O$, let $\vec n_{{\mathcal O}}(x)$ denote the field of unit outward normals to the hypersurfaces $\partial V^t_{\mathcal O }$ and let $s_{\mathcal O}(x) \in \mathbb R$ be the time such that $f^{s_{\mathcal O}(x)}(x)\in V^0_{\mathcal O}$. We define the vector field $\vec v'(x)$ on $\mathbb S^3\times \mathbb S^1$ by the formulas
$$\vec v'(x) = \begin{cases}
	(1 - s_A(x))(1 + s_A(x))\vec n(x) + s_A(x)\vec v(x), & x\in V_A^{-1} \setminus V^1_A \\
	(1 - s_R(x))(1 + s_R(x))\vec n(x) + s_R(x)\vec v(x), & x\in V_R^1 \setminus V^{-1}_R \\
	v(x), & \text{ otherwise}
\end{cases}$$
and denote by $f'^t$ the flow he induced on $\mathbb S^3\times \mathbb S^1$.
For $\delta\in\{-1,+1\}$ we define the diffeomorphism $w_\delta\colon \mathbb S^3\times \mathbb S^1\to \mathbb S^3\times \mathbb S^1 $ by the formula
$$w_\delta(x, y) = (x, \delta y).$$
For $C=(\delta_A, \delta_R)\in \mathbb S^0\times\mathbb S^0$ we define the flow $f^t_C$ on $\mathbb S^3\times \mathbb S^1$ by the formula
$$f^t_{C}(x) = \begin{cases}
	w_{\delta_A} f'^t w_{\delta_A} (x), & x\in V_A\\
	w_{\delta_R} f'^t w_{\delta_R} (x), & x\in V_R\\
	f^t, & \text{ otherwise.}
\end{cases}$$
It is easy to see that the flow's $f^t_{C}$ invariant is $C$ and its ambient manifold is $\mathbb S^3\times \mathbb S^1$.
\end{proof}

\section{Ambient manifolds of the flow of class $G^-_3(M^4)$}
In this section, we prove~\Cref{th:topology}: the only 4-manifold admitting $G^-_3(M^4)$ flows is $\mathbb S^3\times\mathbb S^1$. Moreover, the set $G^-_3(\mathbb S^3\times\mathbb S^1)$ consists of exactly eight equivalence classes of the considered flows.

\begin{proof} Assume that $f^t\in G^-_3(M^4)$ then by~\Cref{th:top-eqv} $f^t$ is topologically equivalent to either the flow $f^t_{C}$ or flow $f^{-t}_{C}$ for some $C\in \mathbb S^0\times\mathbb S^0$. And since the topological equivalence of flows implies the homeomorphism of their ambient manifolds, the supporting manifold of the flow $f^t$ is homeomorphic to $\mathbb S^3\times \mathbb S^1$.
	
	Since the elements of the set $C\in \mathbb S^0\times\mathbb S^0$ correspond one-to-one to the equivalence classes of flows from $G^{-1}_3(M^4)$, the family $G_3^-( M^4) = G^{-2}_3(M^4)\sqcup G^{-1}_3(M^4)$ contains $8$ topological equivalence classes.
\end{proof}

\bibliographystyle{plain} 
\bibliography{refs_En} 

\begin{thebibliography}{10}

\bibitem{Azimov}
Daniel Asimov.
\newblock Round handles and non-singular morse-smale flows.
\newblock {\em Annals of Mathematics}, 102(1):41--54, 1975.

\bibitem{twoLinks}
Christian Bonatti, Viatcheslav Grines, and Elisabeth Pecou.
\newblock Two-dimensional links and diffeomorphisms on 3-manifolds.
\newblock {\em Ergodic Theory and Dynamical Systems}, 22(3):687--710, 2002.

\bibitem{CamposCorderoMartinezAlfaroVindel}
Beatriz Campos, Alicia Cordero, Jose Mart{\'\i}nez~Alfaro, and P~Vindel.
\newblock Nms flows on three-dimensional manifolds with one saddle periodic
  orbit.
\newblock {\em Acta Mathematica Sinica}, 20(1):47--56, 2004.

\bibitem{begin}
Vyacheslav Grines, Timur Medvedev, and Olga Pochinka.
\newblock {\em Dynamical systems on 2-and 3-manifolds}, volume~46.
\newblock Springer, 2016.

\bibitem{Irwin}
M.C. Irwin.
\newblock A classification of elementary cycles.
\newblock {\em Topology}, 9(1):35--47, 1970.

\bibitem{Kosn}
Czes Kosniowski.
\newblock {\em A first course in algebraic topology}.
\newblock Cambridge University Press, 1980.

\bibitem{lickorish}
WB~Raymond Lickorish.
\newblock Homeomorphisms of non-orientable two-manifolds.
\newblock In {\em Mathematical Proceedings of the Cambridge Philosophical
  Society}, volume~59, pages 307--317. Cambridge University Press, 1963.

\bibitem{PoSh-wild}
Olga Pochinka and Danila Shubin.
\newblock On 4-dimensional flows with wildly embedded invariant manifolds of a
  periodic orbit.
\newblock {\em Applied Mathematics and Nonlinear Sciences}, 5(2):261--266,
  2020.

\bibitem{PoSh}
Olga Pochinka and Danila Shubin.
\newblock Non-singular morse--smale flows on n-manifolds with
  attractor--repeller dynamics.
\newblock {\em Nonlinearity}, 35(3):1485, 2022.

\bibitem{PoSh-eqv}
Olga Pochinka and Danila Shubin.
\newblock Nonsingular morse--smale flows with three periodic orbits on
  orientable $3$-manifolds.
\newblock {\em Mathematical Notes}, 112(3):436--450, 2022.

\bibitem{PoSh-top}
Olga Pochinka and Danila Shubin.
\newblock Topology of ambient 3-manifolds of non-singular flows with twisted
  saddle orbit.
\newblock {\em arXiv preprint arXiv:2212.13224}, 2022.

\bibitem{Rolfsen}
Dale Rolfsen.
\newblock {\em Knots and links}, volume 346.
\newblock American Mathematical Soc., 2003.

\bibitem{Shu21}
Danila Shubin.
\newblock Topology of ambient manifolds of nonsingular flows with three twisted
  orbits (in russian).
\newblock {\em Izvestiya Vysshikh uchebnykh zavedeniy. Prikladnaya nelineynaya
  dinamika}, 29(6):863--868, 2021.

\bibitem{Sm}
Stephen Smale.
\newblock Differentiable dynamical systems.
\newblock {\em Bulletin of the American mathematical Society}, 73(6):747--817,
  1967.

\bibitem{Umansky}
Ya.~L. Umanskii.
\newblock Necessary and sufficient conditions for topological equivalence of
  three-dimensional morse--smale dynamical systems with a finite number of
  singular trajectories.
\newblock {\em Matematicheskii Sbornik}, 181(2):212--239, 1990.

\end{thebibliography}
\end{document}